\documentclass[11pt]{article}
\usepackage{}
\usepackage{amssymb}
\usepackage{amsmath}
\usepackage{amsthm}
\usepackage{amsfonts}
\usepackage{mathrsfs}
\newtheorem{theorem}{Theorem}[section]

\newtheorem{lemma}{Lemma}[section]

\newtheorem{definition}{Definition}[section]

\newtheorem{proposition}{Proposition}[section]
\newtheorem{condition}{Condition}[section]

\numberwithin{equation}{section} \allowdisplaybreaks[4]
\begin{document}
\title{Zero-sum continuous-time Markov pure jump game over a fixed duration}
\author{Xin Guo \thanks{Department of Mathematical Sciences, University of
Liverpool, Liverpool, L69 7ZL, U.K.. E-mail: X.Guo21@liv.ac.uk.} and Yi
Zhang \thanks{Corresponding author. Department of Mathematical Sciences, University of
Liverpool, Liverpool, L69 7ZL, U.K.. E-mail: yi.zhang@liv.ac.uk.}}
\date{}
\maketitle

\par\noindent{\bf Abstract:}  This paper considers a two-person zero-sum continuous-time Markov pure jump game in Borel state and action spaces over a fixed finite horizon. The main assumption on the model is the existence of a drift function, which bounds the reward rate. Under some regularity conditions, we show that the game has a value, and both of the players have their optimal policies.
\bigskip

\bigskip


\section{Introduction}
In this paper, we consider a zero-sum continuous-time Markov pure jump game over a fixed finite horizon in Borel state and action spaces.

Zero-sum stochastic games in discrete-time have been intensively studied since 1950s. The first work on this topic is by Shapley, see \cite{Shapley:1953}, where essentially a discounted model with two players in finite state and action spaces was considered, and it was shown that both players have an optimal policy. First extensions of \cite{Shapley:1953} to possibly uncountable state and action spaces are \cite{Maitra:1970,Maitra:1971}, which studied discounted games and positive games, respectively. The case of admissible action spaces dependent on the current state was incorporated in \cite{Parthasarathy:1973}. The state spaces were required to be compact in \cite{Maitra:1970,Maitra:1971,Parthasarathy:1973}. This restriction was relaxed in later works, see e.g., \cite{Couwenbergh:1980,Kumar:1981}. All the aforementioned works consider Borel measurable policies at the cost of extra regularity conditions. The more general class of universally measurable policies in the context of stochastic games was studied intensively in \cite{Nowak:1985}. All these works deal with reward functions bounded or bounded from one side, except \cite{Couwenbergh:1980}, which seems to be the first one considering, in the context of stochastic games, Borel models with a reward function unbounded from both above and below, but bounded by a drift function. The interested reader is referred to \cite{Filar:1997,Neyman:2003} for a comprehensive review of the development in stochastic games in discrete-time.

In continuous-time, one of the first works on stochastic game is \cite{Zachrisson:1964} considering a model in finite state and action spaces. Later, there have been extensions in various directions. Of great relevance to this paper are those on continuous-time Markov pure jump games, i.e., roughly speaking, the process under control is of pure jump type in a Borel state space. Zero-sum continuous-time Markov pure jump game was studied in \cite{GuoGame:2007} with the discounted criterion, and in \cite{GuoGame:2011} with the average criterion.
Both \cite{GuoGame:2007,GuoGame:2011} consider the problem over an infinite time horizon. See also the recent book \cite{PrietoRumeau:2012} for more relevant works on this topic.

The present paper considers a two-person zero-sum continuous-time Markov pure jump game in Borel state and action spaces over a fixed finite time horizon. Our main contributions are the following. We show that there exists a value of the game, and both the maximizer and minimizer have their optimal Markov policies. If the state space is denumerable, similar results were obtained in \cite{Wei:2016}. However, the argument in \cite{Wei:2016}, see e.g., the proof of Theorem 4.1 therein, which is similar to the proof of Theorem 4.1 in \cite{GuoHuangHuang:2015} dealing with a single player case, is essentially based on the denumerable state space, and does not carry to the general Borel model. To deal with the general Borel model, we follow a different and more transparent argument, which is based on the value iteration.  The value iteration argument is also used in \cite{GuoGame:2007} for a discounted model over an infinite horizon. Compared with \cite{GuoGame:2007}, the finite horizon model becomes technically more demanding, because instead of dealing with the space of probability measures as in \cite{GuoGame:2007}, one now needs deal with the space of stochastic kernels, and for that we make use of the Young topology, see the proof of Proposition \ref{GuoXinZhangGameProposition} below. Finally, we mention that in this paper, as in \cite{GuoGame:2007,GuoGame:2011,Wei:2016}, the reward rate is assumed to be bounded by a drift function, but compared with them, our model is nonhomogeneous, i.e., the transition and reward rates are age-dependent, and with a terminal reward.

The rest of this paper is organized as follows. In Section \ref{GuoXinZhangGameSec2} we describe the zero-sum continous-time Markov pure jump game model under consideration. In Section \ref{GuoXinZhangGameSec3}, we impose conditions needed in this paper, and present their immediate consequences. The main statement is presented and proved in Section \ref{GuoXinZhangGameSec4}. This paper is finished with a conclusion in Section \ref{GuoXinZhangGameSec5}.

\section{Model description}\label{GuoXinZhangGameSec2}

In what follows, ${\cal{B}}(X)$ is
the Borel $\sigma$-algebra of the topological space $X,$ $I$ stands for the indicator function, and $\delta_{x}(\cdot)$
is the Dirac measure concentrated on the singleton $\{x\},$ assumed to be measurable. A measure is $\sigma$-additive and $[0,\infty]$-valued. Below, unless stated otherwise, the term of
measurability is always understood in the Borel sense. For a Borel space $X$, let $\mathbb{P}(X)$ be the space of Borel probability measures on $X$. We endow $\mathbb{P}(X)$ with the standard weak topology, and then $\mathbb{P}(X)$ is a Borel space, too. If $X$ is a compact Borel space, then so is $\mathbb{P}(X)$. Throughout
this paper, we adopt the conventions of
$
\frac{0}{0}:=0,~0\cdot\infty:=0,~\frac{1}{0}:=+\infty,~\infty-\infty:=\infty.
$ Finally, the left limit of a function say $f$ at $t\in(-\infty,\infty)$, provided that it exists, is denoted by $f(t-).$

Let $S$ be a nonnempty Borel state space, $A$ be a nonempty Borel action space for the maximizer, $B$ be a nonempty Bore action space for the minimizer. Unless stated otherwise, we endow any Borel space with their Borel $\sigma$-algebras, and below measurability is understood in the Borel sense. Given the current state $x\in S$ and time $t\in [0,\infty)$, the set of admissible actions of the maximizer is $A(t,x)\subseteq A$, and the one of the minimizer is $B(t,x)\subseteq B.$ It is assumed that for each $x\in S$ and $t\in[0,\infty)$, $A(t,x)\subseteq A(x)$ and $B(t,x)\subseteq B(x)$ for some compact subsets $A(x)$ and $B(x)$ of $A$ and $B$, respectively. Furthermore, we assume that the multifunctions $(t,x)\in [0,\infty)\times S\rightarrow A(t,x)$ and $(t,x)\in [0,\infty)\times S\rightarrow B(t,x)$ are both nonempty and compact-valued, and measurable. Consequently, by Theorem 3 of \cite{Himmelberg:1976}, their graphs
$Gr(A)=\{(t,x,a):a\in A(t,x)\}$ and $Gr(B)=\{(t,x,b):b\in B(t,x)\}$ are measurable subsets of $[0,\infty)\times S\times A$ and $[0,\infty)\times S\times B$, respectively. It follows that the set $\mathbb{K}:=\{(t,x,a,b):a\in A(t,x),~b\in B(t,x)\}$ is a measurable subset of $[0,\infty)\times S\times A\times B.$

The transition rate $q(dy|t,x,a,b)$ is a (measurable) signed kernel on $S$ given $(t,x,a,b)\in\mathbb{K}$ such that for each $(t,x,a,b)\in \mathbb{K},$ $q(S|t,x,a,b)=0;$ for each $t\in[0,\infty)$ and $x\in S,$  $q(\Gamma|t,x,a,b)\in[0,\infty)$ for each $\Gamma\in {\cal B}(S)$ such that $x\notin \Gamma;$ and for each $x\in S,$
\begin{eqnarray*}
q_x:=\sup_{t\in[0,\infty),a\in A(t,x),~b\in B(t,x)}q(S\setminus\{x\}|t,x,a,b)<\infty.
\end{eqnarray*}
For convenience, the following notations are used below:
\begin{eqnarray*}
\tilde{q}(dy|t,x,a,b):=q(dy\setminus\{x\}|t,x,a,b),~q(t,x,a,b):=\tilde{q}(S|t,x,a,b),~\forall~(t,x,a,b)\in\mathbb{K}.
\end{eqnarray*}
The reward rate paid by the minimizer to the maximizer at $(t,x,a,b)\in\mathbb{K}$ is $r(t,x,a,b)$, where $r$ is a real-valued measurable function on $\mathbb{K}$. We shall consider the zero-sum game over a fixed time horizon. At the end of the duration say $T>0$, if the state is $x\in S$, then there is a terminal reward of $g(T,x)$ being paid by the minimizer to the maximizer, where $g$ is a real-valued measurable function on $[0,\infty)\times S.$

Roughly speaking, the game is played as follows,  c.f. p.229 of \cite{PrietoRumeau:2012}. At the current time $t$, both players observe the current state $x\in S$, as well as the past states and jump moments of the system. They independently choose some action $a\in A(t,x)$ and $b\in B(t,x)$ according to their policies. Then over a small time increment $[t,t+dt]$, the maximizer receives $r(t,x,a,b)dt$ from the minimizer. The process makes a transition from $x\in S$ to $\Gamma\in {\cal B}(S)$ not consisting of $x$ with probability $q(\Gamma|t,x,a,b)dt +o(dt)$, and the process stays in $x\in S$ over the small time increment with probability $1-q(t,x,a,b)dt+o(dt).$ At the terminal time $T$, if the state of the process is $y\in S,$ then the minimizer pays $g(T,y)$ to the maximizer.

We briefly describe the controlled process and policies of the two players as follows. The more details can be found in Chapter 4 of \cite{Kitaev:1995}. Let us take the sample space $\Omega$ by adjoining to the
countable product space $S\times((0,\infty)\times S)^\infty$ the
sequences of the form
$(x_0,\theta_1,\dots,\theta_n,x_n,\infty,x_\infty,\infty,x_\infty,\dots),$
where $x_0,x_1,\dots,x_n$ belong to $S$,
$\theta_1,\dots,\theta_n$ belong to $(0,\infty),$ and
$x_{\infty}\notin S$ is the isolated point. We equip $\Omega$ with
its Borel $\sigma$-algebra $\cal F$.

Let $t_0(\omega):=0=:\theta_0,$ and for each $n\geq 0$, and each
element $\omega:=(x_0,\theta_1,x_1,\theta_2,\dots)\in \Omega$, let
\begin{eqnarray*}
t_n(\omega)&:=&t_{n-1}(\omega)+\theta_n,
\end{eqnarray*}
and
\begin{eqnarray*}
t_\infty(\omega):=\lim_{n\rightarrow\infty}t_n(\omega).
\end{eqnarray*}
Obviously, $t_n(\omega)$ are measurable mappings on $(\Omega,{\cal
F})$. In what follows, we often omit the argument $\omega\in
\Omega$ from the presentation for simplicity. Also, we regard
$x_n$ and $\theta_{n+1}$ as the coordinate variables, and note
that the pairs $\{t_n,x_n\}$ form a marked point process with the
internal history $\{{\cal F}_t\}_{t\ge 0},$ i.e., the filtration
generated by $\{t_n,x_n\}$. The marked point process $\{t_n,x_n\}$
defines the stochastic process $\{\xi_t,t\ge 0\}$ on $(\Omega,{\cal F})$
 by
\begin{eqnarray}\label{ZhangExponentialGZCTMDPdefxit}
\xi_t=\sum_{n\ge 0}I\{t_n\le t<t_{n+1}\}x_n+I\{t_\infty\le
t\}x_\infty.
\end{eqnarray}
Here we accept $0\cdot x:=0$ and $1\cdot x:=x$ for each $x\in S_\infty,$ and below we denote
$S_{\infty}:=S\bigcup\{x_\infty\}$.

\begin{definition}
\begin{itemize}
\item[(a)]
A policy $\pi$ for the maximizer is given by a
sequence $\{\pi_n\}_{n=0}^\infty$ such that, for each $n=0,1,2,\dots,$
$\pi_n(da|x_0,\theta_1,\dots,x_{n},s)$ is a stochastic kernel on
$A$ concentrated on $A(t_n+s,x_n)$ given $x_0\in S,~\theta_1\in (0,\infty),\dots,~x_{n}\in S,~s\in(0,\infty)$, and for each $\omega=(x_0,\theta_1,x_1,\theta_2,\dots)\in
\Omega$, $t> 0,$
\begin{eqnarray*}
\pi(da|\omega,t)=I\{t\ge t_\infty\}\delta_{a_\infty}(da)+
\sum_{n=0}^\infty I\{t_n< t\le
t_{n+1}\}\pi_{n}(da|x_0,\theta_1,\dots,\theta_n,x_n, t-t_n),
\end{eqnarray*}
where $a_\infty\notin A$ is some isolated point.  A policy $\psi$ for the minimizer is defined as the one for the maximizer, where $\pi$, $A$, $a$ and $da$ are replaced by $\psi$, $B$, $b$ and $db$.
\item[(b)] A policy $\pi$ (or $\psi$) for the maximizer (or the minimizer) is called Markov if one can write $\pi(da|\omega,t)=\pi^M(da|\xi_{t-},t)$ (or $\psi(da|\omega,t)=\psi^M(da|\xi_{t-},t)$) whenever $t<t_\infty$ for some stochastic kernel $\pi^M$ (or $\psi^M$) on $A$ (or $B$) concentrated on $A(t,x)$ (or $B(t,x)$) from $(x,t)\in S\times (0,\infty).$ A Markov policy $\pi$ (or $\psi$) is identified with the underlying stochastic kernel $\pi^M$ (or $\psi^M$).

\end{itemize}
\end{definition}
The class of all policies for the maximizer is denoted by $\Pi,$ and the class of all policies for the minimizer is denoted by $\Psi$.

Under a pair of policies $(\pi,\psi)\in \Pi\times\Psi$, we define the
following random measure on $S\times (0,\infty)$
\begin{eqnarray*}
\nu^{\pi,\psi}(dt, dy)&:=& \int_{A\times B}\tilde{q}(dy|t,\xi_{t-}(\omega),a,b)\pi(da|\omega,t)\psi(db|\omega,t)dt
\end{eqnarray*}
with $q(t,x_\infty,a_\infty,b_\infty):=0=:q(dy|t,x_\infty,a_\infty,b_\infty).$ Then, given the initial distribution  $\gamma$, i.e., a probability measure on ${\cal B}(S)$,
there exists a unique probability measure ${P}^{\pi,\psi}_\gamma$ such
that
\begin{eqnarray*}
{P}_{\gamma}^{\pi,\psi}(x_0\in dx)=\gamma(dx),
\end{eqnarray*}
and with
respect to $P_\gamma^{\pi,\psi},$ $\nu^{\pi,\psi}$ is the dual predictable
projection of the random measure associated with the marked point
process $\{t_n,x_n\}$; see \cite{Jacod:1975,Kitaev:1995}.
When $\gamma$ is a Dirac measure concentrated at $x\in S,$
we use the notation ${}{P}_x^{\pi,\psi}.$ Expectations with respect to
$ {P}_\gamma^{\pi,\psi}$ and ${}{P}_x^{\pi,\psi}$ are denoted as
$ {E}_{\gamma}^{\pi,\psi}$ and ${}{E}_{x}^{\pi,\psi},$ respectively.

The following remark follows from \cite{Jacod:1975}.
Under a pair of policies $(\pi,\psi)\in\Pi\times\Psi$, with the initial distribution $\gamma$, the conditional distribution of $(\theta_{n+1},x_{n+1})$ with the condition on $x_0,\theta_1,x_1,\dots,\theta_{n},x_n$ is given on $\{\omega:x_n(\omega)\in S\}$ by
\begin{eqnarray*}
&&P_\gamma^{\pi,\psi}(\theta_{n+1}\in \Gamma_1,~x_{n+1}\in \Gamma_2|x_0,\theta_1,x_1,\dots,\theta_{n},x_n)\\
&=&\int_{\Gamma_1}e^{-\int_0^t \int_{A\times B} q(s,x_n,a,b)\pi_n(da|x_0,\theta_1,\dots,\theta_n,x_n,s)\psi_n(db|x_0,\theta_1,\dots,\theta_n,x_n,s) ds}\\
&&\int_{A\times B}\tilde{q}(\Gamma_2|t,x_n,a,b)\pi_n(da|x_0,\theta_1,\dots,\theta_n,x_n,t)\psi_n(db|x_0,\theta_1,\dots,\theta_n,x_n,t) dt,\\
&&~\forall~\Gamma_1\in{\cal B}((0,\infty)),~\Gamma_2\in{\cal B}(S);\\
&&P_\gamma^{\pi,\psi}(\theta_{n+1}=\infty,~x_{n+1}=x_\infty|x_0,\theta_1,x_1,\dots,\theta_{n},x_n)\\
&=&e^{-\int_0^\infty  \int_{A\times B} q(s,x_n,a,b)\pi_n(da|x_0,\theta_1,\dots,\theta_n,x_n,s)\psi_n(db|x_0,\theta_1,\dots,\theta_n,x_n,s)ds},
\end{eqnarray*}
and given on $\{\omega:x_n(\omega)=x_\infty\}$ by
\begin{eqnarray*}
P_\gamma^{\pi,\psi}(\theta_{n+1}=\infty,~x_{n+1}=x_\infty|x_0,\theta_1,x_1,\dots,\theta_{n},x_n)=1.
\end{eqnarray*}

Now let $T\in(0,\infty)$ be a fixed time duration, and put
\begin{eqnarray*}
W(x,\pi,\psi):=E_x^{\pi,\psi}\left[\int_0^T \int_{A\times B}r(t,\xi_t,a,b)\pi(da|\omega,t)\psi(db|\omega,t)dt\right]+E_x^{\pi,\psi}\left[g(T,\xi_T)\right]
\end{eqnarray*}
for each $(\pi,\psi)\in\Pi\times\Psi,$ and $x\in S.$ The conditions to be imposed below assure that the above expectations are finite, see Lemma \ref{GuoXinZhangGameLem1}.

The lower value of the zero-sum continuious-time Markov pure jump game over the fixed horizon $[0,T]$ is defined by
\begin{eqnarray*}
L(x):=\sup_{\pi\in\Pi}\inf_{\psi\in\Psi}W(x,\pi,\psi),~\forall~x\in S,
\end{eqnarray*}
and the upper value is defined by
\begin{eqnarray*}
U(x):=\inf_{\psi\in\Psi}\sup_{\pi\in\Pi}W(x,\pi,\psi),~\forall~x\in S.
\end{eqnarray*}
Apparently, $U(x)\ge L(x)$ for each $x\in S.$
If $U(x)=L(x)$ for each $x\in S,$ the function $V$ defined by their common values is called the value of the game.

\begin{definition}
A policy $\pi^\ast\in \Pi$ is called optimal for the maximizer if $\inf_{\psi\in \Psi}W(x,\pi^\ast,\psi)=U(x)$ for each $x\in S.$ A policy $\psi^\ast\in\Psi$ is called optimal for the minimizer if $\sup_{\pi\in \Pi}W(x,\pi,\psi^\ast)=L(x)$ for each $x\in S.$
\end{definition}
It follows that the pair of optimal policies $(\pi^\ast,\psi^\ast)$ in the above definition satisfies
\begin{eqnarray*}
U(x)=\inf_{\psi\in \Psi}W(x,\pi^\ast,\psi)\le W(x,\pi^\ast,\psi^\ast)\le \sup_{\pi\in \Pi}W(x,\pi,\psi^\ast)=L(x),~\forall~x\in S.
\end{eqnarray*}
Then $U(x)=L(x)$ for each $x\in S,$ i.e., the value of the game exists, if both players have their own optimal policies.

The main objective of this paper is to show, under some conditions, that the function $V$ exists, and both players have an optimal policy.

\section{Conditions and relevant facts}\label{GuoXinZhangGameSec3}

In this section, we present the conditions imposed on the continuous-time Markov pure jump game model, and formulate their relevant consequences.

\begin{condition}\label{GuoXinZhangYiCon1}
There exist $[1,\infty)$-valued measurable functions $w_0$ and $w_1$ on $S$ and real constants
$c_0>0$, $c_1>0$, $M_0>0$ and $M_1>0$ such that the following assertions hold.
\begin{itemize}
\item[(a)] For each $(t,x,a,b)\in\mathbb{K}$,  $\int_S w_0(y)q(dy|t,x,a,b)\le c_0w_0(x).$

\item[(b)] For each $x\in S,$ $q_x\le M_0 w_0(x)$.

\item[(c)] For each $(t,x,a,b)\in\mathbb{K}$, $|r(t,x,a,b)|\le M_0 w_0(x)$, $|g(t,x)|\le M_0 w_0(x)$.

\item[(d)] For each $(t,x,a,b)\in\mathbb{K}$, $\int_S w_1(y)q(dy|t,x,a,b)\le c_1w_1(x).$

\item[(e)] For each $x\in S,$ $w_0(x)q_x\le M_1 w_1(x)$.
 \end{itemize}
\end{condition}

Consider a fixed $[1,\infty)$-valued measurable function say $f$ on $S$. A
function $u$ on $[0,T]\times S$ is called $f$-bounded if
\begin{eqnarray*}
\|u\|_{f}:=\sup_{(t,x)\in [0,T]\times S}\frac{|u(t,x)|}{f(x)}<\infty.
\end{eqnarray*}
The set of $f$-bounded measurable functions on $[0,T]\times S$ is denoted by $B_{f}([0,T]\times S).$ Fix a function $u\in B_{f}([0,T]\times S).$ Suppose for each $x\in S,$ $u(\cdot,x)$ is absolutely continuous. Then it is known that there exists a measurable function $u'$ on $[0,T]\times S$ such that for each $x\in S,$ the derivative of $u(x,\cdot)$ exists and coincides with $u'(\cdot,x)$ almost everywhere on $[0,T]$.

Under Condition \ref{GuoXinZhangYiCon1}, let $C_{w_0,w_1}^{1,0}([0,T]\times S)$ be the collection of functions $u\in
B_{w_0}([0,T]\times S)$ such that for each $x \in S$, $u(\cdot,x)$ is
absolutely continuous, and $u'$ belongs to
$B_{w_0+w_1}([0,T]\times S)$.

\begin{lemma}\label{GuoXinZhangGameLem1}
Suppose Condition \ref{GuoXinZhangYiCon1} is satisfied. Let some pair of policies $(\pi,\psi)\in\Pi\times \Psi$ be arbitrarily fixed. Then the following assertions hold.
\begin{itemize}
\item[(a)] $P_{x}^{\pi,\psi}(t_\infty=\infty)=1$ for each $x\in S.$
\item [(b)] $E_{x}^{\pi,\psi}[w_0(\xi_{t})] \le e^{c_0t}w_0(x)$ for each $t\ge 0$ and $x\in S$.
\item[(c)] $ |W(x,\pi,\psi)|\le (T+1) M_0e^{c_0T}w_0(x) $ for each $x\in S$.
\item[(d)] For each $u\in C_{w_0,w_1}^{1,0}([0,T]\times S)$,
 \begin{eqnarray*}
&& E_{x}^{\pi,\psi}\left[\int_{0}^{T}
\left(u'(t,\xi_t)+\int_S \int_{A}\int_B u(t,x)q(dx|t,\xi_t,a,b)\pi(da|\omega,t)\psi(db|\omega,t)\right)dt
\right]\\
&=&E_{x}^{\pi,\psi}[u(T,\xi_{T})]-u(0,x).
\end{eqnarray*}
for each $x\in S.$
\end{itemize}
\end{lemma}
\par\noindent\textit{Proof.} See Lemmas 3.1, 3.2 and 3.3 in \cite{GuoHuangZhang:2016}. $\hfill\Box$
\bigskip

Throughout the rest of this paper, let $m$ be an $[1,\infty)$-valued measurable function on $S$ such that $q_x\le m(x)$ for each $x\in S.$ Such a function exists by the Novikov seperation theorem, see \cite{Kechris:1995}. We introduce the following stochastic kernel on $S$ from $(t,x,a,b)\in \mathbb{K}$ defined by
\begin{eqnarray*}
\tilde{p}(dy|t,x,a,b):=\delta_x(dy)+\frac{q(dy|t,x,a,b)}{m(x)},~\forall~(t,x,a,b)\in \mathbb{K}.
\end{eqnarray*}

 \begin{condition}\label{GuoXinZhang:2016Con2} For each $t\in [0,T]$ and $x\in S$,
    \begin{itemize}
\item[(a)]   $r(t,x,a,b)$ is  continuous in $(a,b)\in A(t,x)\times B(t,x)$; and

\item[(b)]  for each measurable function $u$ on $S$ such that $\sup_{x\in S}\frac{|u(x)|}{w_0(x)}<\infty,$ $\int_S u(y)\tilde{p}(dy|t,x,a,b)$ is continuous in $(a,b)\in A(t,x) \times B(t,x)$.
  \end{itemize}
  \end{condition}

Suppose that Condition \ref{GuoXinZhangYiCon1} is satisfied. For each $t\in[0,\infty)$, $x\in S$, $\lambda\in  \mathbb{P}(A(t,x))$ and $\mu\in\mathbb{P}(B(t,x))$, we introduce the notations
\begin{eqnarray*}
q(dy|t,x,\lambda,\mu)&:=&\int_{A(t,x)}\int_{B(t,x)} q(dy|t,x,a,b)\lambda(da)\mu(db), \\
r(t,x,\lambda,\mu)&:=&\int_{A(t,x)}\int_{B(t,x)}r(t,x,a,b)\lambda(da)\mu(db).
\end{eqnarray*}
(In particular, the integral in the second line of the above is finite under Condition \ref{GuoXinZhangYiCon1}.)  Then 
$q(dy|t,x,\lambda,\mu)$ and
$r(t,x,\lambda,\mu)$ are measurable on ${\cal K},$ where
\begin{eqnarray*}
{\cal K}:=\left\{(t,x,\lambda,\mu)\in [0,\infty)\times S\times \mathbb{P}(A)\times \mathbb{P}(B): \lambda\in \mathbb{P}(A(t,x)),~\mu\in \mathbb{P}(B(t,x))\right\}.
\end{eqnarray*}
In greater details, since $(t,x)\rightarrow A(t,x)$ and $(t,x)\rightarrow B(t,x)$ are measurable and compact-valued multifunctions, as assumed earlier, by Theorem 3 of \cite{Himmelberg:1975} and Proposition 7.22 of \cite{Bertsekas:1978}, so are the multifunctions $(t,x)\rightarrow \mathbb{P}(A(t,x))$ and $(t,x)\rightarrow \mathbb{P}(B(t,x)).$ It follows from Theorem 3 of \cite{Himmelberg:1976} that
 ${\cal K}$ is measurable in the Borel space $[0,\infty)\times S\times \mathbb{P}(A)\times \mathbb{P}(B)$.
By Corollary 7.29.1 and Lemma 7.21 of \cite{Bertsekas:1978} that $q(dy|t,x,\lambda,\mu)$ and
$r(t,x,\lambda,\mu)$ are measurable on ${\cal K}.$

The next lemma, used repeatedly in the next section, is known. But we include its rather short proof for completeness. Recall that $A(t,x)$ and $B(t,x)$ are compact subsets of $A$ and $B$ as assumed in the beginning of the model description.
\begin{lemma}\label{GuoXinZhangGameLem2}  Suppose that Conditions \ref{GuoXinZhangYiCon1} and \ref{GuoXinZhang:2016Con2} are satisfied.
\begin{itemize}
\item[(a)]  Let $t\in[0,T]$ and $x\in S$ be arbitrarily fixed. For each  $u\in B_{w_0}([0,T]\times S)$,  the functions $r(t,x,\lambda,\mu)$ and $\int_S u(t,y)q(dy|t,x,\lambda,\mu)$ are continuous in $(\lambda,\mu)\in \mathbb{P}(A(t,x))\times \mathbb{P}(B(t,x))$.

\item[(b)]  If a function $h(t,x,\lambda,\mu)$ is real-valued and measurable on ${\cal K},$ and continuous in $(\lambda,\mu) \in \mathbb{P}(A(t,x))\times \mathbb{P}(B(t,x))$ (for each fixed $(t,x)\in [0,T]\times S$), then the function
\begin{eqnarray*}
(t,x,\lambda)\rightarrow \inf_{\mu \in \mathbb{P}(B(t,x))} h(t,x,\lambda,\mu)
\end{eqnarray*}
is measurable on $\{(t,x,\lambda)\in[0,T]\times S\times \mathbb{P}(A):\lambda\in \mathbb{P}(A(t,x))\}$ and continuous in $\lambda \in \mathbb{P}(A(t,x))$ (for each fixed $(t,x)\in [0,T]\times S$).
\end{itemize}
\end{lemma}
\par\noindent\textit{Proof.} (a) For the fixed $t\in[0,T]$ and $x\in S$, the functions $r(t,x,a,b)$ and $\int_S u(t,y)q(dy|t,x,a,b)$ are bounded and continuous in $(a,b)\in A(t,x)\times  B(t,x)$. The statement follows from Corollary 7.29.1 and Lemma 7.12 of \cite{Bertsekas:1978}, and the Tietze extension theorem.

(b) The first assertion follows from Theorem 2 of \cite{Himmelberg:1976}.  The second assertion is a consequence of the Berge theorem, see Theorem 17.31 in \cite{Aliprantis:2006}. $\hfill\Box$
\bigskip

\section{Main statement}\label{GuoXinZhangGameSec4}
In this section, we present and prove the main result of this paper; see Theorem \ref{GuoXinZhangGameTheorem} below.

Under Conditions \ref{GuoXinZhangYiCon1} and \ref{GuoXinZhang:2016Con2}, it follows from Lemmas \ref{GuoXinZhangGameLem1} and \ref{GuoXinZhangGameLem2} and the fundamental theorem of calculus that
the following operator $G$ maps $u\in B_{w_0}([0,T]\times S)$ to $C_{w_0,w_1}^{1,0}([0,T]\times S)$:
\begin{eqnarray*}
 &&G[u](t,x)\\
 &:=& e^{-m(x)(T-t)}g(T,x)+\int_{0}^{T-t}e^{-m(x)s}\\
 &&\sup_{ \lambda\in \mathbb{P}(A(t+s,x))}\inf_{\mu\in \mathbb{P}(B(t+s,x))}\left\{r(t+s,x,\lambda,\mu)
 +m(x)\int_S u(t+s,y)\tilde{p}(dy|t+s,x,\lambda,\mu)\right\}ds
\end{eqnarray*}
for each $t\in[0,T]$ and $x\in S.$

\begin{proposition}\label{GuoXinZhangGameProposition}
Suppose that Conditions \ref{GuoXinZhangYiCon1} and \ref{GuoXinZhang:2016Con2} are satisfied.
There is a fixed point of the operator $G$ in $C_{w_0,w_1}^{1,0}([0,T]\times S)$.
\end{proposition}
\par\noindent\textit{Proof.} Let us define
\begin{eqnarray*}
u_0(t,x):=\frac{M_0}{c_0}\left\{c_0e^{c_0(T-t)}+e^{c_0(T-t)}-1\right\} w_0(x)\ge 0
\end{eqnarray*}
for each $t\in [0,T]$ and $x\in S.$ Then $u_0$ belongs to $C_{w_0,w_1}^{1,0}([0,T]\times S).$ For each $n\ge 0$, we legitimately define $u_{n+1}:=G[u_n]$. The rest of the proof goes in two steps.

\underline{\textit{Step 1.}} Show that $\{u_n\}$ is a monotone nonincreasing sequence, and for each $n=0,1,\dots,$ 
\begin{eqnarray*}
|u_n(t,x)|\le u_0(t,x)=\frac{M_0}{c_0}\left\{c_0e^{c_0(T-t)}+e^{c_0(T-t)}-1\right\} w_0(x).
\end{eqnarray*}
for each $t\in[0,T]$ and $x\in S.$

For each $t\in[0,T]$ and $x\in S,$
\begin{eqnarray*}
&&u_1(t,x)=G[u_0](t,x)\\
&\le& e^{-m(x) (T-t)}M_0w_0(x)+\int_0^{T-t} e^{-m(x)s}\\
&&\sup_{ \lambda\in \mathbb{P}(A(t+s,x))}\inf_{\mu\in \mathbb{P}(B(t+s,x))}\left\{M_0w_0(x)
 +m(x)\int_S u_0(t+s,y)\tilde{p}(dy|t+s,x,\lambda,\mu)\right\}ds\\
 &=&e^{-m(x) (T-t)}M_0w_0(x)+M_0w_0(x)\int_0^{T-t} e^{-m(x)s}ds\\
 &&+\frac{M_0}{c_0}\int_0^{T-t} e^{-m(x)s} m(x)\left\{c_0e^{c_0(T-t-s)}+e^{c_0(T-t-s)}-1\right\} w_0(x)ds \\
&& +\frac{M_0}{c_0}\int_0^{T-t} e^{-m(x)s}\sup_{ \lambda\in \mathbb{P}(A(t+s,x))}\inf_{\mu\in \mathbb{P}(B(t+s,x))}\left\{(c_0e^{c_0(T-t-s)}+e^{c_0(T-t-s)}-1)\right.\nonumber\\
&&  \left.\int_S w_0(y)q(dy|t+s,x,\lambda,\mu)\right\}ds   \\
&\le&  e^{-m(x) (T-t)}M_0w_0(x)+M_0w_0(x)\int_0^{T-t} e^{-m(x)s}ds\\
 &&+\frac{M_0}{c_0}\int_0^{T-t} e^{-m(x)s} m(x)\left\{c_0e^{c_0(T-t-s)}+e^{c_0(T-t-s)}-1\right\} w_0(x)ds \\
&& +\frac{M_0}{c_0}\int_0^{T-t} e^{-m(x)s}(c_0e^{c_0(T-t-s)}+e^{c_0(T-t-s)}-1) c_0 w_0(x)ds,
\end{eqnarray*}
where the first and the  last inequalities are by Condition \ref{GuoXinZhangYiCon1}.
For the third summand on the right hand side of the last inequality, integration by parts gives
\begin{eqnarray*}
&&\frac{M_0}{c_0}\int_0^{T-t} e^{-m(x)s} m(x)\left\{c_0e^{c_0(T-t-s)}+e^{c_0(T-t-s)}-1\right\} w_0(x)ds\\
&=&-w_0(x)M_0 e^{-m(x)(T-t)}+\frac{M_0}{c_0}\left\{c_0 e^{c_0(T-t)}+e^{c_0(T-t)}-1\right\}w_0(x)\\
&&-\frac{M_0}{c_0}w_0(x)\int_0^{T-t} e^{-m(x)s}\left\{c_0^2 e^{c_0(T-t-s)}+c_0e^{c_0(T-t-s)}\right\}ds.
\end{eqnarray*}
This, together with the previous calculations, shows that
\begin{eqnarray*}
u_1(t,x)\le \frac{M_0}{c_0}\left\{c_0 e^{c_0(T-t)}+e^{c_0(T-t)}-1\right\}w_0(x)=u_0(t,x),~\forall~t\in[0,T],~x\in S.
\end{eqnarray*}
It follows from this and the monotonicity of the operator $G$ that $\{u_n\}$ is a monotone nonincreasing sequence, and for each $n\ge 0,$
\begin{eqnarray*}
u_n(t,x)\le u_0(t,x)=\frac{M_0}{c_0}\left\{c_0e^{c_0(T-t)}+e^{c_0(T-t)}-1\right\} w_0(x)
\end{eqnarray*}
for each $t\in[0,T]$ and $x\in S.$

On the other hand, a similar calculation to the above gives
\begin{eqnarray*}
&&u_1(t,x)\\
&\ge& -e^{-m(x) (T-t)}M_0w_0(x)-\int_0^{T-t} e^{-m(x)s}\\
&&\sup_{ \lambda\in \mathbb{P}(A(t+s,x))}\inf_{\mu\in \mathbb{P}(B(t+s,x))}\left\{-M_0w_0(x)
 -m(x)\int_S u_0(t+s,y)\tilde{p}(dy|t+s,x,\lambda,\mu)\right\}ds\\
 &\ge &-u_0(t,x)
\end{eqnarray*}
for each $t\in[0,T]$ and $x\in S.$ Hence, for each $n\ge 0$,
\begin{eqnarray*}
|u_n(t,x)|\le u_0(t,x)=\frac{M_0}{c_0}\left\{c_0e^{c_0(T-t)}+e^{c_0(T-t)}-1\right\} w_0(x).
\end{eqnarray*}
for each $t\in[0,T]$ and $x\in S.$

\underline{\textit{Step 2.}} Consider the function $u^\ast$ defined by $u^\ast(t,x):=\lim_{n\rightarrow\infty}u_n(t,x)$ for each $t\in [0,T]$ and $x\in S$. The limit exists due to monotone convergence. We show that $u^\ast$ is a fixed point of the operator $G$ in $C_{w_0,w_1}^{1,0}([0,T]\times S)$.

It follows from the definition of $u^\ast$ and what was established in Step 1 that
\begin{eqnarray*}
|u^\ast(t,x)|\le \frac{M_0}{c_0}\left\{c_0e^{c_0(T-t)}+e^{c_0(T-t)}-1\right\} w_0(x)
\end{eqnarray*}
for each $t\in[0,T]$ and $x\in S,$ that is, $u^\ast\in B_{w_0}([0,T]\times S).$

We verify that $u^\ast$ is a fixed point of $G$ as follows. It is evident that for each $n\ge 0,$ $G[u^\ast](t,x)\le G[u_n](t,x)=u_{n+1}(t,x)$ for each $t\in[0,T]$ and $x\in S.$ Hence,
\begin{eqnarray}\label{GuoXinZhangGame:06}
G[u^\ast](t,x) \le u^\ast(t,x)
\end{eqnarray}
for each $t\in[0,T]$ and $x\in S.$

The rest of this proof mainly verifies the opposite direction of the above inequality. Let $x\in S$ be fixed, and consider the space of $\mathbb{P}(A)$-valued measurable mappings say $\lambda$ on $[0,T]$ such that for each $t \in [0,T]$, $\lambda_t\in \mathbb{P}(A(t,x)).$ We denote this space by ${\cal P}_A(x)$. The notation ${\cal P}_B(x)$ is understood similarly, with $A$ being replaced by $B.$

Note that by Theorem 2 of \cite{Himmelberg:1976}, applicable due to Lemma \ref{GuoXinZhangGameLem2}, for each $x\in  S$ and $t\in[0,T]$,
\begin{eqnarray*}
&&\int_{0}^{T-t}e^{-m(x)s}\\
 &&\sup_{ \lambda\in \mathbb{P}(A(t+s,x))}\inf_{\mu\in \mathbb{P}(B(t+s,x))}\left\{r(t+s,x,\lambda,\mu)
 +m(x)\int_S u(t+s,y)\tilde{p}(dy|t+s,x,\lambda,\mu)\right\}ds\\
&=&\int_{0}^{T-t}e^{-m(x)s}\\
 &&\sup_{ \lambda\in {\cal P}_A(x)}\inf_{\mu\in {\cal P}_B(x)}\left\{r(t+s,x,\lambda_{t+s},\mu_{t+s})
 +m(x)\int_S u(t+s,y)\tilde{p}(dy|t+s,x,\lambda_{t+s},\mu_{t+s})\right\}ds.
\end{eqnarray*}

Fix $(t,x)\in[0,T]\times S$ and some $\mu\in {\cal P}_B(x)$ arbitrarily. By Theorem 2 of \cite{Himmelberg:1976} and Lemma \ref{GuoXinZhangGameLem2}, for each $n\ge 0$, there exists $\lambda^n\in {\cal P}_A(x)$ such that
\begin{eqnarray}\label{GuoXinZhangGame:02}
&&u_{n+1}(t,x)=G[u_n](t,x)\nonumber\\
&\le& e^{-m(x)(T-t)}g(T,x)+\int_{0}^{T-t}e^{-m(x)s}\nonumber\\
 && \left\{r(t+s,x,\lambda^n_{t+s},\mu_{t+s})
 +m(x)\int_S u_n(t+s,y)\tilde{p}(dy|t+s,x,\lambda^n_{t+s},\mu_{t+s})\right\} ds
\end{eqnarray}

Now we endow the space ${\cal P}_A(x)$ with the Young topology. The interested reader is referred to Section 44 of \cite{Davis:1993} or Chapter 2 of \cite{Costa:2013} for the details; or p.249-250 of \cite{BauerleRieder:2011} for a brief description. The relevant fact to this paper is that endowed with this topology, the space ${\cal P}_A(x)$ is compact metrizable, and the mapping
\begin{eqnarray*}
\lambda\in {\cal P}_A(x)\rightarrow \int_0^T \int_{A(x)} g(t,a)\lambda_{t}(da)dt
\end{eqnarray*}
is continuous for each measurable real-valued function $g$ on $[0,T]\times A(x)$ such that $t\in[0,T]\rightarrow g(t,a)$ is measurable,  $a\in A(x)\rightarrow g(t,a)$ is continuous, and $\int_0^T \sup_{a\in A(x)}|g(t,a)|dt<\infty.$  Recall that $A(t,x)\subseteq A(x)$ for each $x\in S$ and $t\in[0,\infty).$

Since ${\cal P}_A(x)$ is compact metrizable, without loss of generality we assume that the sequence $\{\lambda^n\}$ in ${\cal P}_A(x)$ converges to some $\lambda^\ast\in {\cal P}_A(x)$, for otherwise one can take a convergent subsequence and relabel it. Note that
\begin{eqnarray}\label{GuoXinZhangGame:01}
&& \left|\int_{0}^{T-t}e^{-m(x)s}
 m(x)\int_S u_n(t+s,y)\tilde{p}(dy|t+s,x,\lambda^n_{t+s},\mu_{t+s}) ds\right.\nonumber\\
 &&\left.-
 \int_{0}^{T-t}e^{-m(x)s}
 m(x)\int_S u^\ast(t+s,y)\tilde{p}(dy|t+s,x,\lambda^n_{t+s},\mu_{t+s}) ds\right|\nonumber\\
 &\le&\int_{0}^{T-t}e^{-m(x)s}
 m(x)\int_{A(t+s,x)}\int_S |u_n(t+s,y)-u^\ast(t+s,y)|\tilde{p}(dy|t+s,x,a,\mu_{t+s})\lambda^n_{t+s}(da) ds\nonumber\\
 &\le&\int_{0}^{T-t}e^{-m(x)s}
 m(x)\sup_{a\in A(t+s,x)}\left\{\int_S |u_n(t+s,y)-u^\ast(t+s,y)|\tilde{p}(dy|t+s,x,a,\mu_{t+s})\right\}ds.\nonumber\\
\end{eqnarray}
On the other hand,
\begin{eqnarray*}
&&\lim_{n\rightarrow\infty}\sup_{a\in A(t+s,x)}\left\{\int_S |u_n(t+s,y)-u^\ast(t+s,y)|\tilde{p}(dy|t+s,x,a,\mu_{t+s})\right\}\\
&=&\sup_{a\in A(t+s,x)}\left\{\lim_{n\rightarrow\infty}\int_S |u_n(t+s,y)-u^\ast(t+s,y)|\tilde{p}(dy|t+s,x,a,\mu_{t+s})\right\}\\
&=&0,
\end{eqnarray*}
where the first equality is by Theorem A.1.5 of \cite{BauerleRieder:2011}, applicable under Condition \ref{GuoXinZhang:2016Con2}, and the last equality is by the dominated convergence theorem, applicable under Condition \ref{GuoXinZhangYiCon1}. It follows from this, (\ref{GuoXinZhangGame:01}) and the dominated convergence theorem that
\begin{eqnarray*}
&&\lim_{n\rightarrow \infty}\left|\int_{0}^{T-t}e^{-m(x)s}
 m(x)\int_S u_n(t+s,y)\tilde{p}(dy|t+s,x,\lambda^n_{t+s},\mu_{t+s}) ds\right.\\
 &&\left.-
 \int_{0}^{T-t}e^{-m(x)s}
 m(x)\int_S u^\ast(t+s,y)\tilde{p}(dy|t+s,x,\lambda^n_{t+s},\mu_{t+s}) ds\right|=0.
\end{eqnarray*}

Now as $n\rightarrow\infty$,
\begin{eqnarray*}
&&\left|\int_{0}^{T-t}e^{-m(x)s}
 m(x)\int_S u_n(t+s,y)\tilde{p}(dy|t+s,x,\lambda^n_{t+s},\mu_{t+s}) ds\right.\\
 &&\left.-
 \int_{0}^{T-t}e^{-m(x)s}
 m(x)\int_S u^\ast(t+s,y)\tilde{p}(dy|t+s,x,\lambda^\ast_{t+s},\mu_{t+s}) ds\right|\\
&\le& \left|\int_{0}^{T-t}e^{-m(x)s}
 m(x)\int_S u_n(t+s,y)\tilde{p}(dy|t+s,x,\lambda^n_{t+s},\mu_{t+s}) ds\right. \\
 &&\left.-
 \int_{0}^{T-t}e^{-m(x)s}
 m(x)\int_S u^\ast(t+s,y)\tilde{p}(dy|t+s,x,\lambda^n_{t+s},\mu_{t+s}) ds\right| \\
 &&+ \left|\int_{0}^{T-t}e^{-m(x)s}
 m(x)\int_S u^\ast(t+s,y)\tilde{p}(dy|t+s,x,\lambda^n_{t+s},\mu_{t+s}) ds\right. \\
 &&\left.-
 \int_{0}^{T-t}e^{-m(x)s}
 m(x)\int_S u^\ast(t+s,y)\tilde{p}(dy|t+s,x,\lambda^\ast_{t+s},\mu_{t+s}) ds\right|\\
&\rightarrow& 0,
\end{eqnarray*}
where the convergence to zero is also by the definition of the Young topology. It follows from this and the definition of the Young topology again that, after passing to the limit as $n\rightarrow \infty$ on the both sides of (\ref{GuoXinZhangGame:02}),
\begin{eqnarray}\label{GuoXinZhangGame:05}
&&u^\ast(t,x)\nonumber\\
&\le& e^{-m(x)(T-t)}g(T,x)+\int_{0}^{T-t}e^{-m(x)s}\left\{ r(t+s,x,\lambda^\ast_{t+s},\mu_{t+s})\right.\nonumber\\
&&\left. +m(x)\int_S u^\ast(t+s,y)\tilde{p}(dy|t+s,x,\lambda^\ast_{t+s},\mu_{t+s})\right\} ds\nonumber\\
&\le& e^{-m(x)(T-t)}g(T,x)+\int_{0}^{T-t}e^{-m(x)s} \sup_{\lambda\in \mathbb{P}(A(t+s,x))}\left\{r(t+s,x,\lambda,\mu_{t+s})\right.\nonumber\\
&&\left. +m(x)\int_S u^\ast(t+s,y)\tilde{p}(dy|t+s,x,\lambda,\mu_{t+s})\right\} ds.
\end{eqnarray}

By Theorem 2 of \cite{Himmelberg:1976}, applicable due to Lemma \ref{GuoXinZhangGameLem2}, there exists $\mu^\ast\in {\cal P}_B(x)$ such that
\begin{eqnarray*}
&&\inf_{\mu\in\mathbb{P}(B(t+s,x))}\sup_{\lambda\in \mathbb{P}(A(t+s,x))}\left\{r(t+s,x,\lambda,\mu)+m(x)\int_S u^\ast(t+s,y)\tilde{p}(dy|t+s,x,\lambda,\mu)\right\}\\
&=&\sup_{\lambda\in \mathbb{P}(A(t+s,x))}\left\{r(t+s,x,\lambda,\mu^\ast_{t+s})+m(x)\int_S u^\ast(t+s,y)\tilde{p}(dy|t+s,x,\lambda,\mu^\ast_{t+s})\right\}
\end{eqnarray*}
for each $s\in[0,T-t].$
By the Ky Fan minimax theorem, see Theorem 2 of \cite{Fan:1953},
\begin{eqnarray*}
&&\sup_{\lambda\in \mathbb{P}(A(t+s,x))}\inf_{\mu\in\mathbb{P}(B(t+s,x))}\left\{r(t+s,x,\lambda,\mu)+m(x)\int_S u^\ast(t+s,y)\tilde{p}(dy|t+s,x,\lambda,\mu)\right\}\\
&=&\sup_{\lambda\in \mathbb{P}(A(t+s,x))}\left\{r(t+s,x,\lambda,\mu^\ast_{t+s})+m(x)\int_S u^\ast(t+s,y)\tilde{p}(dy|t+s,x,\lambda,\mu^\ast_{t+s})\right\}
\end{eqnarray*}
for each $s\in[0,T-t].$
Since $\mu\in {\cal P}_B(x)$ in (\ref{GuoXinZhangGame:05}) was arbitrarily fixed, we see from (\ref{GuoXinZhangGame:05}) and the previous equality that
\begin{eqnarray*}
&&u^\ast(t,x)\\
&\le&  e^{-m(x)(T-t)}g(T,x)+\int_{0}^{T-t}e^{-m(x)s} \sup_{\lambda\in \mathbb{P}(t+s,x)}\left\{r(t+s,x,\lambda,\mu^\ast_{t+s})\right.\nonumber\\
&&\left. +m(x)\int_S u^\ast(t+s,y)\tilde{p}(dy|t+s,x,\lambda,\mu^\ast_{t+s})\right\} ds\\
&=& e^{-m(x)(T-t)}g(T,x)+\int_{0}^{T-t}e^{-m(x)s}\\
&& \sup_{\lambda\in \mathbb{P}(A(t+s,x))}\inf_{\mu\in\mathbb{P}(B(t+s,x))}\left\{r(t+s,x,\lambda,\mu)+m(x)\int_S u^\ast(t+s,y)\tilde{p}(dy|t+s,x,\lambda,\mu)\right\} ds\\
&=&G[u^\ast](t,x).
\end{eqnarray*}
Since $(t,x)\in[0,T]\times S$ was arbitrarily fixed, this and (\ref{GuoXinZhangGame:06}) imply
\begin{eqnarray*}
u^\ast(t,x)=G[u^\ast](x,t),~\forall~t\in[0,T],~x\in S.
\end{eqnarray*}

Finally, since $u^\ast\in B_{w_0}([0,T]\times S)$, and $G$ maps each element of $B_{w_0}([0,T]\times S)$ to $C_{w_0,w_1}^{1,0}([0,T]\times S)$ as mentioned earlier, it follows that $u^\ast$ is a fixed point of $G$ in $C_{w_0,w_1}^{1,0}([0,T]\times S).$
$\hfill\Box$
\bigskip

\begin{theorem}\label{GuoXinZhangGameTheorem}
Suppose that Conditions \ref{GuoXinZhangYiCon1} and \ref{GuoXinZhang:2016Con2} are satisfied.  Then the zero-sum continuous-time Markov pure jump game has a value $V$, and both the maximizer and minimizer have an optimal Markov policy. In particular, there is a pair of Markov policies $(\pi_*^M,\psi_*^M)\in \Pi\times \Psi$ such that $W(x,\pi_*^M,\psi_*^M)=V(x)$ for each $x\in S.$
\end{theorem}

\par\noindent\textit{Proof.}  By Proposition \ref{GuoXinZhangGameProposition}, we can consider a solution $u\in C_{w_0,w_1}^{1,0}([0,T]\times S)$ to the following equation
\begin{eqnarray*}\label{GuoXinZhangGame:07}
&&u(t,x)\nonumber\\
&=&e^{-m(x)(T-t)}g(T,x)+\int_{0}^{T-t}e^{-m(x)s}\nonumber\\
 &&\sup_{ \lambda\in \mathbb{P}(A(t+s,x))}\inf_{\mu\in \mathbb{P}(B(t+s,x))}\left\{r(t+s,x,\lambda,\mu)
 +m(x)\int_S u(t+s,y)\tilde{p}(dy|t+s,x,\lambda,\mu)\right\}ds,\nonumber\\
 &&\forall~t\in [0,T],~x\in S.
\end{eqnarray*}
Then
\begin{eqnarray*}
&&e^{-m(x)t}u(t,x)\\
&=&e^{-m(x)T}g(T,x)+\int_{0}^{T-t}e^{-m(x)(t+s)}\\
 &&\sup_{ \lambda\in \mathbb{P}(A(t+s,x))}\inf_{\mu\in \mathbb{P}(B(t+s,x))}\left\{r(t+s,x,\lambda,\mu)
 +m(x)\int_S u(t+s,y)\tilde{p}(dy|t+s,x,\lambda,\mu)\right\}ds\\
&=&e^{-m(x)T}g(T,x)+\int_{t}^{T}e^{-m(x)(s)}\\
 &&\sup_{ \lambda\in \mathbb{P}(A(s,x))}\inf_{\mu\in \mathbb{P}(B(s,x))}\left\{r(s,x,\lambda,\mu)
 +m(x)\int_S u(s,y)\tilde{p}(dy|s,x,\lambda,\mu)\right\}ds\\
 &&\forall~t\in [0,T],~x\in S.
\end{eqnarray*}
It follows that for each $x\in S,$
\begin{eqnarray}\label{GuoXinZhangGame:08}
u(T,x)=g(T,x)
\end{eqnarray}
and
\begin{eqnarray*}
&&u'(t,x)+\sup_{ \lambda\in \mathbb{P}(A(t,x))}\inf_{\mu\in \mathbb{P}(B(t,x))}\left\{r(t,x,\lambda,\mu)
 +\int_S u(s,y)q(dy|t,x,\lambda,\mu)\right\}=0
\end{eqnarray*}
almost everywhere on $[0,T].$

By Theorem 2 of \cite{Himmelberg:1976}, applicable due to Lemma \ref{GuoXinZhangGameLem2}, there exists a Markov policy say $\pi^M_*$ for the maximizer such that for each $x\in S,$
\begin{eqnarray*}
&&u'(t,x)+\inf_{\mu\in \mathbb{P}(B(t,x))}\left\{\int_A r(t,x,a,\mu)\pi^M_*(da|x,t)
 +\int_S u(t,y)\int_A q(dy|t,x,a,\mu)\pi^M_*(da|x,t)\right\}=0
\end{eqnarray*}
almost everywhere on $[0,T],$ that is, for each $\mu\in \mathbb{P}(B(t,x)),$
\begin{eqnarray*}
&&u'(t,x)+ \int_A r(t,x,a,\mu)\pi^M_*(da|x,t)
 +\int_S u(t,y)\int_A q(dy|t,x,a,\mu)\pi^M_*(da|x,t) \ge 0
\end{eqnarray*}
almost everywhere on $[0,T].$

Now, by Lemma \ref{GuoXinZhangGameLem1}(d), for  each policy $\psi \in \Psi$ for the minimizer and $x\in S,$
\begin{eqnarray*}
&&E_x^{\pi_*^M,\psi}[g(T,\xi_T)]-u(0,x)=E_x^{\pi_*^M,\psi}[u(T,\xi_T)]-u(0,x)\\
&=&E_{x}^{\pi_*^M,\psi}\left[\int_{0}^{T}
\left(u'(t,\xi_t)+\int_S \int_{A}\int_B u(t,x)q(dx|t,\xi_t,a,b)\pi_*^M(da|\xi_t,t)\psi(db|\omega,t)\right)dt
\right]\\
&\ge&-E_{x}^{\pi_*^M,\psi}\left[\int_{0}^{T} \int_{A}\int_B r(t,\xi_t,a,b)\pi_*^M(da|\xi_t,t)\psi(db|\omega,t)dt
\right],
\end{eqnarray*}
where the first equality is by (\ref{GuoXinZhangGame:08}). That is,
\begin{eqnarray*}
u(0,x)\le W(x,\pi^M_*,\psi),~\forall~x\in S.
\end{eqnarray*}
Since $\psi\in \Psi$ was arbitrarily fixed, we see
\begin{eqnarray}\label{GuoXinZhangGame:09}
u(0,x)\le \inf_{\psi\in\Psi}W(x,\pi^M_*,\psi)\le \sup_{\pi\in\Pi}\inf_{\psi\in\Psi}W(x,\pi,\psi)=L(x),~\forall~x\in S.
\end{eqnarray}

Similarly, by By Theorem 2 of \cite{Himmelberg:1976} and the Ky Fan minimax theorem (see Theorem 2 of \cite{Fan:1953}), there exists a Markov policy say $\psi^M_*$ for the minimizer such that for each $x\in S,$
\begin{eqnarray*}
&&u'(t,x)+\sup_{\lambda\in \mathbb{P}(A(t,x))}\left\{\int_B r(t,x,\lambda,b)\psi^M_*(db|x,t)
 +\int_S u(t,y)\int_B q(dy|t,x,\lambda,b)\psi^M_*(db|x,t)\right\}=0
\end{eqnarray*}
almost everywhere on $[0,T].$ Then by using Lemma \ref{GuoXinZhangGameLem1}(d), one can show as in the above that
\begin{eqnarray*}
u(0,x)\ge \sup_{\pi\in\Pi}W(x,\pi,\psi^M_*)\ge \inf_{\psi\in\Psi} \sup_{\pi\in\Pi}W(x,\pi,\psi)=U(x),~\forall~x\in S.
\end{eqnarray*}
Combining this and (\ref{GuoXinZhangGame:09}) yields
\begin{eqnarray*}
u(0,x)=L(x)=U(x)= \sup_{\pi\in\Pi}W(x,\pi,\psi^M_*)=\inf_{\psi\in\Psi}W(x,\pi^M_*,\psi)=W(x,\pi^M_*,\psi_*^M),~\forall~x\in S.
\end{eqnarray*}
The proof is completed.
 $\hfill\Box$
\bigskip

\section{Conclusion}\label{GuoXinZhangGameSec5}
In this paper, we consider a zero-sum continuous-time Markov pure jump game over a fixed finite horizon. The state and action spaces are Borel. Under some conditions, we establish the existence of the value of the game, and both players have a Markov optimal policy. Finally we mention that here we imposed the strong continuity condition, and showed that the value function is measurable. But the same argument applies to show that the value function is continuous (in the state), if we modify the continuity condition more or less obviously and use the corresponding versions of measurable selection theorem in the argument. We omitted the details to avoid repetition.

\section*{Acknowledgement} This work was partially carried out with a financial grant from the Research Fund for Coal and Steel of
the European Commission, within the INDUSE-2-SAFETY project (Grant No. RFSR-CT-2014-00025).

\end{document}